# DPG* METHOD

Brendan Keith[1], Leszek Demkowicz[1] and Jay Gopalakrishnan[2]

[1]Institute for Computational Engineering and Sciences (ICES), The University of Texas at Austin
[2] Fariborz Maseeh Department of Mathematics + Statistics, Portland State University

**Abstract**

We introduce a cousin of the DPG method - the DPG* method - discuss their relationship and compare the two methods through numerical experiments.

## 1 Introduction

The *Ideal Discontinuous Petrov Galerkin (DPG) Method with Optimal Test Functions* admits three interpretations [8]. It can be viewed as a Petrov-Galerkin (PG) discretization in which the (optimal) test functions are computed on the fly. The word *optimal* refers to the fact that the functions realize the supremum in the discrete inf-sup stability condition and, therefore, the PG discretization automatically inherits the stability of the continuous method. The DPG method can also be interpreted as a *Minimum Residual Method* in which the residual is measured in a dual norm implied by the test norm. Finally, the DPG method can also be viewed as a *mixed method* where one solves simultaneously for the Riesz representation of the residual, the so-called *error representation function*. All three formulations involve inversion of the Riesz operator on the (Hilbert) test space which, in general, cannot be done exactly and has to be approximated. This has led to the introduction of a ("larger") *enriched* or *search* test space and an approximate inverse of the Riesz operator obtained using the standard Bubnov-Galerkin method [1]. The corresponding *practical DPG method* retains its three interpretations although we deal now with *approximate* optimal test functions, residuals[2] and error respresentation function.

We mention that the word "discontinuous" in DPG method corresponds to the use of discontinuous test functions (broken test spaces) which makes the whole methodology computationally efficient. The broken test spaces formulation provides a foundation for the DPG methodology, and can be developed for any well posed variational formulation [4].

In the presented philosophy, the mixed problem formulation occurs only for the ideal DPG method and, after the enriched space approximation, for the practical DPG method, but not for the original problem. A different philosophy has been used by Cohen, Dahmen and Welpert in [6]. The original problem

$$Bu = l \tag{1.1}$$

---

[1]An optimal choice for the discretization of an inner product.
[2]The supremum in the dual norm is computed over the large but finite-dimensional test space.

is replaced *at the continuous level* with the mixed problem,

$$\begin{cases} R_V \psi + Bu &= l \\ B'\psi &= 0 \,. \end{cases} \quad (1.2)$$

Here $U, V$ and Hilbert trial and test spaces, $B : U \to V'$ is the operator corresponding to a bilinear (sesquilinear) form $b(u,v)$, $u \in U$, $v \in V$, defining the variational problem, and linear (antilinear) $l \in V'$ is the load. Finally, $R_V : V \to V'$ is the Riesz operator corresponding to the test inner product $(v, \delta v)_V$. Operator $B$ is bounded below,

$$\|Bu\|_{V'} \geq \gamma \|u\|_U, \quad u \in U\,,$$

and form $l$ satisfies the usual compatibility condition: $l \in \mathcal{N}(B')^\perp$, i.e.

$$l(v) = 0 \quad \forall v \in \mathcal{N}(B')$$

where $B' : V'' \sim V \to U'$ is the conjugate operator of $B$. By the Closed Range Theorem, problem (1.1) is well-posed. Function $\psi \in V$ is the Riesz representation of zero residual $l - Bu$ and it is equal zero. The mixed problem (1.2) is trivially equivalent to the original problem. We mention that by operator $B$ we may mean the broken version of the original variational formulation. This means that the original variational formulation has been replaced with

$$\begin{cases} u \in U, \hat{u} \in \hat{U} \\ b(u,v) + \langle \hat{u}, v \rangle = l(v), \quad v \in V(\mathcal{T}) \end{cases} \quad (1.3)$$

Here $V(\mathcal{T})$ is a broken test space corresponding to a Finite Element (FE) partition, we silently assume that the original bilinear form can be extended to a broken test space, $\hat{u}$ represents the additional unknown (Lagrange multiplier) living in an appropriate *trace space* $\hat{U}$, and $\langle \hat{u}, v \rangle$ stands for the coupling term defined over the mesh skeleton, see [4] for details. Introducing the group variable $(u, \hat{u})$ and replacing $b(u,v)$ with the modified form $b_{\text{mod}}((u, \hat{u}), v) := b(u,v) + \langle \hat{u}, v \rangle$, we extend the validity of our discussion to the broken test space setting. The important point of the broken test space setting is that the test inner product must be well defined on the broken test space (a *localizable inner product*).

The error representation function can be eliminated leading to a normal equation for $u$:

$$B' R_V^{-1} Bu = B' R_V^{-1} l \,. \quad (1.4)$$

For broken test spaces, the elimination *can be done locally at the element level*.

Instead of discretizing the original problem (1.1), we discretize the equivalent mixed problem (1.2). The key point is that the discrete test space need not have same dimension as the discrete trial space, *it can be bigger*. In other words, we make sense how to solve an *overdetermined discrete problem* in which $\dim V_h > \dim U_h$. The stability of the mixed discretization is now governed by the two Brezzi inf-sup conditions. The *inf-sup in kernel condition* is satisfied trivially, and the Babuška-Brezzi (BB) condition is *easier to satisfy* since we can work with a larger discrete test space. If we use broken test spaces, the additional cost of forming and solving the discrete equivalent of normal problem (1.4) is local. The discrete BB condition is established through the construction of local [9, 4, 11] or global [5] Fortin operators.

Once we have established the discrete stability of the mixed problem, we have the standard convergence result,

$$\underbrace{\|\psi_h\|_V + \|u - u_h\|_U}_{\text{approximation error}} \leq C(\gamma_h) \underbrace{\inf_{w_h \in U_h} \|u - w_h\|_U}_{\text{best approximation error}} \quad (1.5)$$

where the global stability constant $C$ depends upon continuity constant $\|B\| = \|b\|$ and discrete BB constant $\gamma_h$. Notice that because $\psi = 0$, the best approximation of $\psi$ is zero,

$$\inf_{\phi_h \in V_h} \|\phi - \phi_h\|_V = 0\,.$$

**DPG$^*$ method.** Once the legitimacy of the mixed problem (1.2) has been established, we can put a load into the second equation,

$$\begin{cases} R_V \psi + Bu &= 0 \\ B'\psi &= g \end{cases} \quad (1.6)$$

where $g \in U'$. This version of the mixed problem represents a strategy for solving the adjoint problem:

$$B'\psi = g\,. \quad (1.7)$$

If the original mixed formulation helps to understand how to cope with an *overdetermined* problem, this formulation clearly establishes means for dealing with an *undetermined problem* as, for the discrete version of it, $\dim V_h > \dim U_h$.

We have arrived at the DPG$^*$ method.

**Remark 1** Obviously, if we are interested in using the DPG$^*$ method to approximate the original problem, we should swap in (1.6) operator $B'$ with $B$. ∎

## 2 The DPG$^*$ Method

We begin with a couple of observations about the new formulation at the continuous level.

- Load $g \in U'$ can be arbitrary. Indeed $\mathcal{N}((B')') = \mathcal{N}(B) = \{0\}$ and, therefore, the compatibility condition $g \perp \mathcal{N}((B')')$ is satisfied trivially.

- The null space of $B'$ may be non-trivial, so equation (1.7) may have multiple solutions. However, the first equation in the mixed formulation (1.6) automatically selects the solution in the orthogonal component $\mathcal{N}(B')^\perp$. Indeed, let

$$\psi = \psi_0 + \psi_1, \qquad \psi_0 \in \mathcal{N}(B'), \quad \psi_1 \in \mathcal{N}(B')^\perp\,.$$

Substituting $\psi$ into the first equation, and testing with $\psi_0$, we get $\|\psi_0\|_V^2 = 0$. The same observation holds on the discrete level, except that the null space $V_0$ is now the null space of the discrete approximation of $B'$ [3].

Now, both original and the new mixed problem differ only in the load. The discrete problem,

$$\begin{cases} \psi_h \in V_h,\ u_h \in U_h \\ (\psi_h, \phi_h)_V + b(u_h, \phi_h) = l(\phi_h) & \phi_h \in V_h \\ b^*(\phi_h, w_h) = g(w_h) & w_h \in U_h \end{cases} \quad (2.8)$$

where $b^*(\psi, w) = \overline{b(w, \psi)}$, covers both DPG and DPG$^*$ methods, we enjoy the same discrete stability and a-priori error estimate,

$$(\|\psi - \psi_h\|_V + \|u - u_h\|_U) \leq C(\gamma_h) \left( \inf_{\phi_h \in V_h} \|\psi - \phi_h\|_V + \inf_{w_h \in U_h} \|u - w_h\|_U \right). \quad (2.9)$$

The main difference between the DPG method and the general case (2.9) is that $\psi$, in general, is no longer zero. For the dual problem: $l = 0, g \neq 0$, $\psi$ represents the discrete solution of the dual problem. For broken test spaces, problem (2.9) is solved using the same technology as for the primary DPG problem. Function $\psi_h$ is statically condensed out element-wise, and we solve the global problem for the Lagrange multiplier $u_h$, determining the ultimate solution $\psi_h$ in the back-substitution step.

Note that "over-enriching" space $V_h$ harms efficiency but *not* discrete stability. Indeed, the ultimate convergence rates are determined now by the best approximation rates for both unknowns: $\psi$ and $u$. In the DPG setting $u$ stands for the group variable including the original unknown (the field variable) and a trace variable resulting from breaking the test space. For strong variational formulations, the trace variable is absent, and the method reduces to the well known FOSLL$^\star$ method [2].

Before we continue with more of a theoretical discussion, we shall focus on a simple numerical example comparing the DPG and DPG$^*$ methods.

### Example: Time-Harmonic Linear Acoustics

Let $\Omega = (0,1)^2 \subset \mathbb{R}^2$ be a unit square domain. We look for a pressure $p$ and velocity $u = (u_1, u_2)$, satisfying the following equations in $\Omega$:

$$\begin{cases} i\omega p + \text{div}\, u = 0 \\ i\omega u + \nabla p = 0, \end{cases} \quad (2.10)$$

and impedance Boundary Condition (BC) on boundary $\Gamma = \partial \Omega$,

$$p - u \cdot n = g$$

---

[3] Analogous situation happens for the stabilized formulation of time-harmonic Maxwell equations [7].

| dp | DPG $L^2$-error | DPG* $L^2$-error | DPG* graph norm error |
|---|---|---|---|
| 0 | 40.57 | 31.57 | 284.28 |
| 1 | 33.77 | 17.03 | 77.65 |
| 2 | 33.51 | 18.50 | 33.04 |
| 3 | 36.44 | 34.58 | 39.32 |
| 4 | 37.20 | 39.47 | 42.17 |
| 5 | 37.32 | 40.45 | 42.78 |
| 6 | 37.38 | 40.72 | 42.96 |

Table 1: Plane wave example. Comparison of DPG and DPG* methods.

where $n$ denotes the outward normal unit vector to $\Gamma$.

We shall focus on the DPG method for the (broken) ultraweak variational formulation. Recall that this version of the DPG method suffers from least pollution (dispersion) error and, in the end of an adaptive process, delivers practically $L^2$-projection. The energy spaces and forms are defined as follows.

$$\begin{aligned}
\mathfrak{u} &:= (p, u, \hat{p}, \hat{u} \cdot n) \in U := (L^2(\Omega) \times (L^2(\Omega))^2 \times H^{1/2}(\Gamma_h) \times H^{-1/2}(\Gamma_h^0) \\
\mathfrak{v} &:= (q, v) \in V := H^1(\mathcal{T}) \times H(\text{div}, \mathcal{T}) \\
b(\mathfrak{u}, \mathfrak{v}) &:= (p, i\omega q + \text{div}_h v) + (u, i\omega v + \nabla_h q) + \langle \hat{u} \cdot n, q \rangle_{\Gamma_h^0} + \langle \hat{p}, v \cdot n \rangle_{\Gamma_h} + \langle p, q \rangle_\Gamma \\
l(\mathfrak{v}) &:= \langle g, q \rangle_\Gamma
\end{aligned} \quad (2.11)$$

We shall compute with the adjoint graph test norm,

$$\begin{aligned}
\|\mathfrak{v}\|_V^2 &:= \|A^*(q, v)\|^2 + \|(q, v)\|^2 \\
&= \|i\omega q + \text{div}_h v\|^2 + \|i\omega v + \nabla_h q\|^2 + \|q\|^2 + \|v\|^2
\end{aligned} \quad (2.12)$$

where $\|\cdot\|$ denotes the $L^2$-norm and $A$ represents the first order system operator. Note that that $A$ is (formally) skew-adjoint, $A^* = -A$.

We used a plane wave propagating at a $40^o$ degrees angle for a manufactured solution with two wavelengths ($\omega = 4\pi$), and a mesh of $2 \times 2$ elements with $p = 3$. We use the exact sequence of the first type which means that the $L^2$-variables $p, u$ are discretized with quadratic polynomials, trace $\hat{p}$ is discretized with continuous cubic polynomials defined on the mesh skeleton $\Gamma_h$, and normal flux $\hat{u}_n$ is discretized with discontinuous quadratics defined on the interior mesh skeleton $\Gamma_h^0$. The discontinuous test functions are discretized with polynomials of order $p + dp = 3 + dp$ where $dp = 0, 1, \ldots, 6$. The results are presented in Table 1. The first column displays the $L^2$-error for the field variables for the DPG method and the third column shows the corresponding numbers for the DPG* method with the error measured in the test norm, i.e. the broken adjoint graph norm. For comparison, the second column displays the error for the DPG* method measured in the $L^2$-norm. We report *relative errors*, i.e. they are computed in percent of the corresponding norms of the solution.

Recall that discrete stability depends only upon the BB constant $\beta_h$,

$$\beta_h = \inf_{\|u_h\|_U=1} \sup_{\|v_h\|_V=1} |b(u_h, v_h)|.$$

With the trial space fixed and test space increasing with growing $dp$, the BB constant can only grow. The best approximation error involves three variables: fields $(u, p)$, traces $\hat{u}, \hat{p}$, and test function $(q, v)$. With fixed $p$ and growing $dp$, the best approximation for $(q, v)$ is decreasing but for the trial variables remains constant. Consequently, the best approximation error in trial variables eventually dominates and we expect the approximation error to reach some asymptotics. This is indeed observed in the first column for the DPG method. The sweet spot is for $dp = 2$ and the error then is slightly growing reaching an asymptotic value of 37 percent. The growth of the error can be explained by the fact that we are looking only at one component of the solution without reporting the error in traces and dual variable. With the stability improving and the best approximation error, the combined error in all variables should diminish but we cannot exclude some redistribution within all variables.

The same comments apply to the DPG∗ method and the error reported in the third column. The sweet spot for the graph norm is again for $dp = 2$ but the error grows only slightly with $dp$ reaching an asymptotic value of 43%.

The most interesting is the behavior of the DPG$^*$ error measured in the $L^2$ norm reported in the second column. The best value of 17% is reached for low value of dp=1. It grows slightly for dp=2 but starting with $dp = 3$, the error doubles and then grows to an asymptotic value of 40%. At first, the behavior seems to be inconsistent with the theory as we know it.

The explanation comes from the behavior of the boundedness below constant for the first order operator $A^*$. At the continuous level,

$$\|A^*(q, v)\| \geq \alpha \|(q, v)\| \quad (q, v) \in H^1(\Omega) \times H(\text{div}, \Omega), q + v \cdot n = 0 \text{ on } \Gamma.$$

Obviously, the inequality holds for any *conforming* discretization of $(q, v)$. In the DPG methodology, the conformity and BC are replaced with *weak conformity*. i.e. orthogonality of jumps to traces $\hat{p}, \hat{u} \cdot n$, i.e.,

$$\sum_{e \subset \Gamma} \int_e (q + v \cdot n) p = 0 \quad \forall p$$
$$\sum_{e \subset \Omega} \int_e [v \cdot n] p = 0 \quad \forall p$$
$$\sum_{e \subset \Omega} \int_e [q] u \cdot n = 0 \quad \forall u.$$

With the polynomial degree of $(p, u)$ fixed and the polynomial degree of $(q, v)$ growing, the corresponding discrete boundedness below constant $\alpha_h$ converges to zero. As,

$$(1 + \alpha_h)\|(q, v)\|^2 \leq \|(q, v)\|^2 + \|A^*(q, v)\|^2,$$

a good approximation in the adjoint graph norm does not longer translate into a good approximation in the $L^2$-norm.

The difference is the $L^2$-norm quality is quite dramatic and important from the practical point of view as, for high wave regime, we are interested mostly in the $L^2$-norm quality of the solution. Fig. 1 presents the DPG and DPG∗ solutions for $dp = 1$. The difference in quality of the two solutions appears in higher wave regime as well although is less dramatic. For a comparison, Fig. 2 provides an analogous comparison for the same problem with 10 wavelengths and a $10 \times 10$ mesh with $p = 3$ and $dp = 2$. Due to pollution, both solutions are underresolved but the DPG∗ solution is still better than the DPG solution. The corresponding $L^2$ errors are 46% for the DPG method and 36% for the DPG∗ method. For $dp = 1$ that yielded the best result for the two wavelengths problem, the difference in quality is even less pronounced (58% vs. 53%). All observations concerning the dependence of the $L^2$-norm error on $dp$ remain valid for the 10-wavelength example. For instance, for $dp = 3$, the $L^2$ error jumps up to 68%.

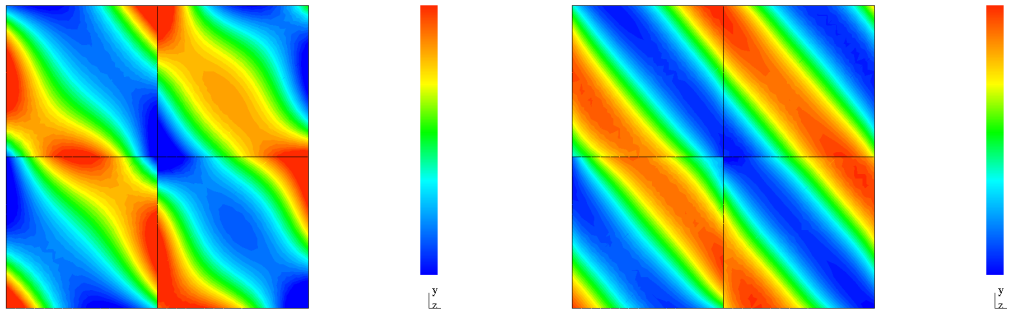

Figure 1: Plane wave example. Comparison of DPG (left) and DPG*(right) solutions (real part of pressure) for the adjoint graph test norm and $dp = 1$.

All presented solutions have been displayed on the same scale, from -1 to 1 and with the same resolution. The under- and overshoots in the DPG solutions are bigger than for its competitor. The difference in uder-resolution patterns visible in Fig.2 can be attributed to the change in sign in the impedance BC. Finally, we emphasize that the cost of obtaining both solutions is identical as the two methods share the same stiffness matrix.

**Convergence rates.** We proceed by checking asymptotic $h$-convergence rates. We use the same plane wave problem for the manufactured solution but with just one wavelength only. We have solved the problem with elements of order $p = 1, 2, 3, 4$ performing four consecutive global $h$-refinements. For each polynomial order, we solve the problem four times, using $dp = 0, 1, 2, 3$. Fig. 3 presents the convergence rates with respect to number of (trial) degrees-of-freedom, obtained using the adjoint graph test norm. The optimal convergence rates deteriorate for $p = 3$ and are lost for $p = 4$ except for the case of $dp = 0$ indicating

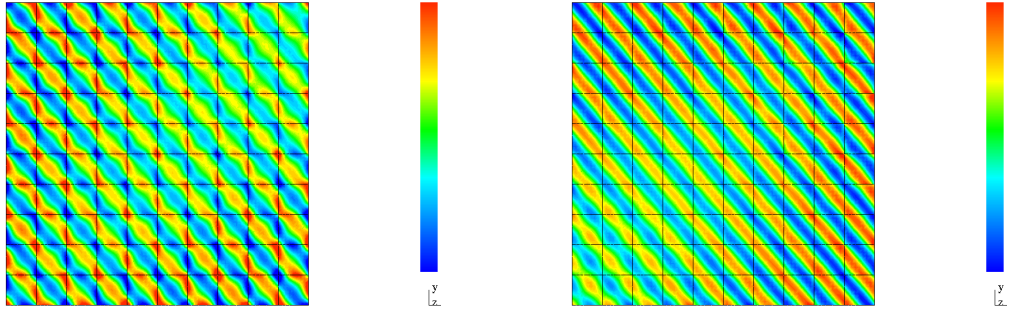

Figure 2: Plane wave example with 10 wavelengths. Comparison of DPG (left) and DPG*(right) solutions for the adjoint graph test norm and $dp = 2$.

conditioning problems. In general, results for $dp = 0$ are shifted upward indicating a smaller stability constant. Fig. 4 presents the same experiment but with the "mathematician's test norm',

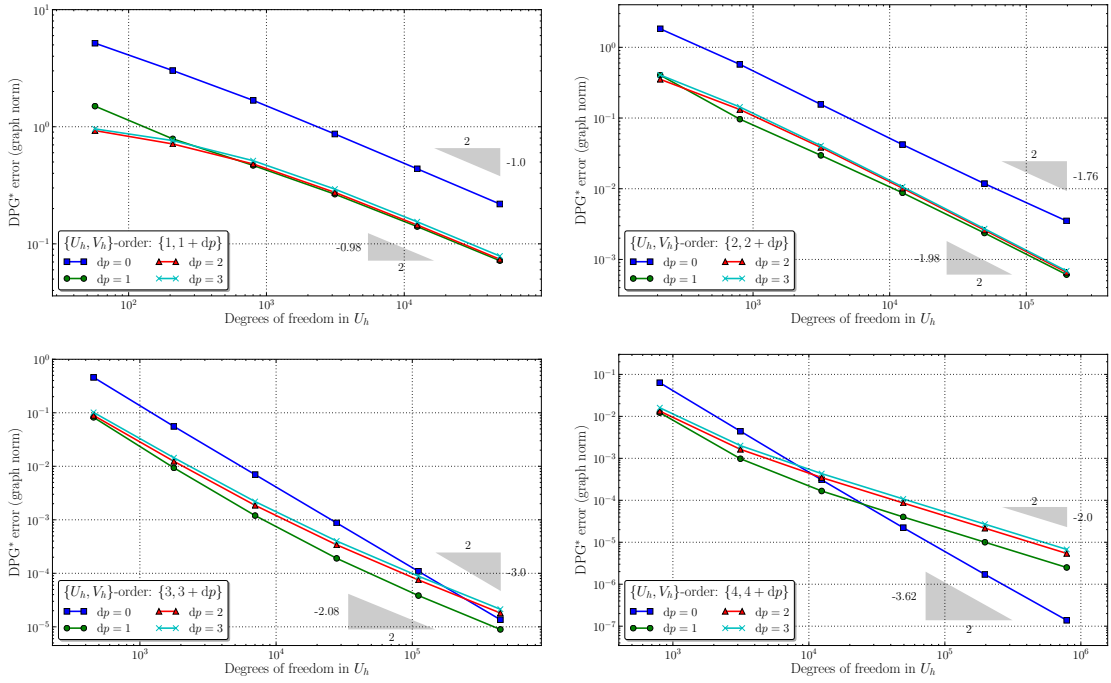

Figure 3: Plane wave example with one wavelength. Adjoint graph test norm. $h$-convergence rates for $p = 1$ (top left), $p = 2$ (top right), $p = 3$ (bottom left), $p = 4$ (bottom right), with $dp = 0, 1, 2, 3$.

$$\|(q, v)\|_V^2 := \|q\|_{H^1(\mathcal{T}_h)}^2 + \|v\|_{H(\mathrm{div}, \mathcal{T}_h)}^2.$$

The results are more or less the same. In the lowest order case $p = 1, dp = 0$, the method takes longer to get into the asymptotic range, and the optimal rates are lost for $p = 3$ and higher $dp$.

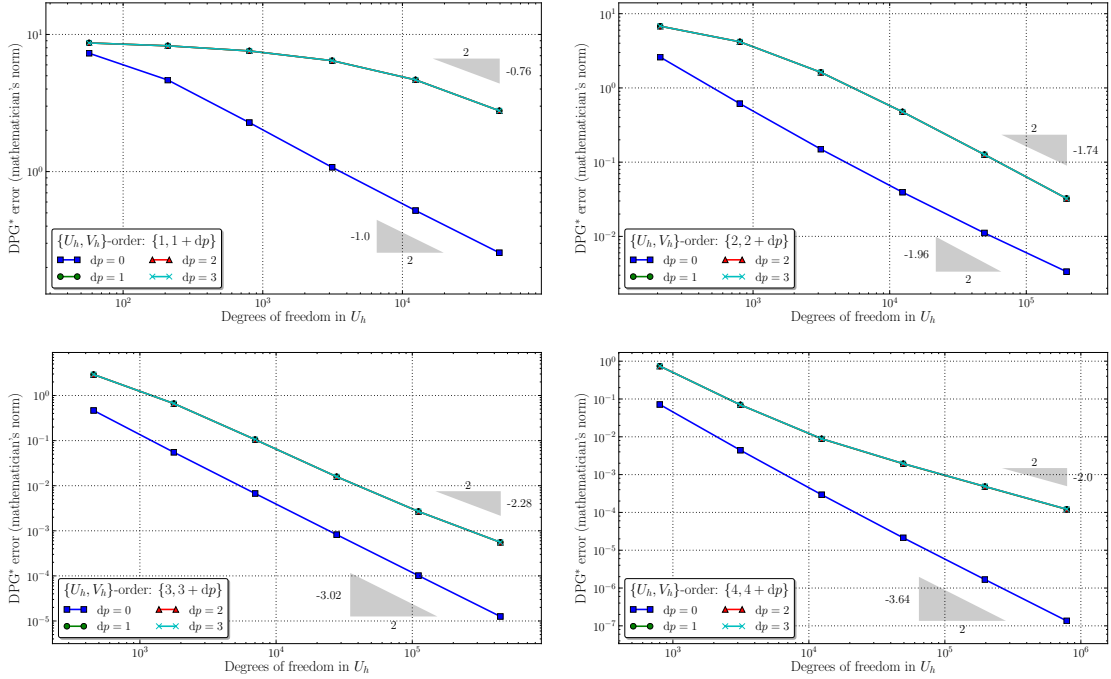

Figure 4: Plane wave example with one wavelength. Mathematician's test norm. $h$-convergence rates for $p = 1$ (top left), $p = 2$ (top right), $p = 3$ (bottom left), $p = 4$ (bottom right), with $dp = 0, 1, 2, 3$.

## 3 Stability Analysis for the Mixed Problem

We have come across the DPG* method when studying a-posteriori goal-oriented estimation and adaptivity. As both original and dual problems are just special cases of the mixed method, it is important to understand the stability of the method. Of course, the general Brezzi's theory applies [1] but the problem of interest and the choice of norms that are used here is non-standard, so the general theory does not yield optimal estimates. Results presented in this section can be found in [12] although, we believe, we use simpler and more "compact" arguments.

Let $U, V$ be two Hilbert spaces, and $B : U \to V'$ be a linear, continuous and *bounded below* operator. We do not assume that $B$ is surjective. Consider the mixed problem:

$$\begin{cases} \psi \in V, \, u \in U \\ R_V \psi + Bu &= l \\ B' \psi &= g \end{cases} \quad (3.13)$$

where $l \in V'$, $g \in U'$ are arbitrary elements of the dual space of antilinear functionals defined on $V$ or $U$, respectively. Operator $R_V : V \to V'$ denotes the Riesz operator for space $V$.

Notice that the mixed problem has a unique solution for any $l$ and $g$. If $l \in \mathcal{R}(B)$ and $g = 0$ then $\psi$ represents the Riesz representation of the (zero) residual, and it is equal to zero as well.

Let $V_0 := \mathcal{N}(B')$. Consider the orthogonal decompositon,

$$V = V_0 \oplus V^\perp, \qquad \psi = \psi_0 + \psi^\perp. \tag{3.14}$$

Notice that $R_V \psi_0$ and $Bu$ are orthogonal in $V'$. Indeed,

$$(Bu, R_V \psi_0)_{V'} = \langle Bu, \psi_0 \rangle = \langle B' \psi_0, u \rangle = 0.$$

Recalling that the Riesz operator is an isometry, the first equation in (3.13),

$$R_V \psi_0 + Bu = l - R_V \psi^\perp,$$

and the Pythagoras Theorem imply that

$$\|\psi_0\|_V^2 + \|Bu\|_{V'}^2 = \|l - R_V \psi^\perp\|_{V'}^2. \tag{3.15}$$

In what follows, we replace the original norm in $U$ with an equivalent "energy norm",

$$\|u\|_E := \|Bu\|_{V'},$$

and equip the dual space $U'$ with the corresponding dual norm. We get,

$$\begin{aligned}
\|g\|_{U'} &= \sup_u \frac{|g(u)|}{\|Bu\|_{V'}} = \sup_u \frac{|\langle B'\psi^\perp, u\rangle|}{\|Bu\|_{V'}} = \sup_u \frac{|\langle Bu, \psi^\perp\rangle|}{\|Bu\|_{V'}} \\
&= \sup_{l \in \mathcal{R}(B)} \frac{|\langle l, \psi^\perp\rangle|}{\|l\|_{V'}} = \sup_{l \in \mathcal{N}(B')^\perp} \frac{|\langle l, \psi^\perp\rangle|}{\|l\|_{V'}} \\
&= \|\psi^\perp\|_V.
\end{aligned} \tag{3.16}$$

Adding (3.16) to (3.15) and using orthogonality of $\psi_0$ and $\psi^\perp$, we obtain the fundamental identity:

$$\|\psi\|_V^2 + \|Bu\|_{V'}^2 = \|l - R_V \psi^\perp\|_{V'}^2 + \|g\|_{U'}^2. \tag{3.17}$$

Without any extra information about the relation of $l$ to $R_V \psi^\perp$, we can use only the triangle inequality to obtain the bound:

$$\|\psi\|_V^2 + \|Bu\|_{V'}^2 \leq (\|l\|_{V'} + \|g\|_{U'})^2 + \|g\|_{U'}^2. \tag{3.18}$$

The stability bound translates immediately into an a-posteriori error estimate for the mixed problem. If $(\psi, u)$ denotes the solution and $(\psi_h, u_h)$ an arbitrary element of $(V \times U)$ then,

$$\|\psi - \psi_h\|_V^2 + \|B(u - u_h)\|_{V'}^2 \leq (\|l - R_V \psi_h - Bu_h\|_{V'} + \|g - B'\psi_h\|_{U'})^2 + \|g - B'\psi_h\|_{U'}^2. \tag{3.19}$$

Estimate (3.19) bounds the error in both $\psi$ and $u$ by the residuals corresponding to the two equations and may serve as a starting point for more practical a-posteriori estimates based on a specific way of estimating the residuals.

If we are interested in estimating only one component, we can obtain slightly better bounds than (3.18). Identity (3.17) implies immediately,

$$\|\psi\|_V^2 + \|Bu\|_{V'}^2 \leq \|l\|_{V'}^2 + 2\|l\|_{V'}\underbrace{\|\psi^\perp\|_V}_{=\|g\|_{U'}} + \|\psi^\perp\|_V^2 + \|g\|_{U'}^2,$$

which, given $\|\psi^\perp\|_V \leq \|\psi\|_V$, immediately implies

$$\|Bu\|_{V'} \leq \|l\|_{V'} + \|g\|_{U'}. \tag{3.20}$$

Again, this translates into the a-posteriori error estimate:

$$\|B(u - u_h)\|_{V'} \leq \|l - R_V\psi_h - Bu_h\|_{V'} + \|g - B'\psi_h\|_{U'}. \tag{3.21}$$

To obtain an estimate for $\psi$ alone, we start with the "cosine law",

$$\|\psi\|_V^2 + \|Bu\|_{V'}^2 + 2\Re\underbrace{\langle \psi, Bu\rangle}_{=\langle \psi^\perp, Bu\rangle} = \|l\|_{V'}^2$$

and use Cauchy-Schwarz and Young's inequalities to obtain,

$$\|\psi\|_V^2 + \|Bu\|_{V'}^2 \leq \|l\|_{V'}^2 + \underbrace{\|\psi^\perp\|_V^2}_{=\|g\|_{U'}^2} + \|Bu\|_{V'}^2$$

and the final estimate,

$$\|\psi\|_V^2 \leq \|l\|_{V'}^2 + \|g\|_{U'}^2. \tag{3.22}$$

Again, this implies the a-posteriori error estimate for error $\psi - \psi_h$,

$$\|\psi - \psi_h\|_V^2 \leq \|l - R_V\psi_h - Bu_h\|_{V'}^2 + \|g - B'\psi_h\|_{U'}^2. \tag{3.23}$$

**A posteriori error estimation for the DPG and the DPG$^*$ methods.** Assuming $g = 0$, we obtain a-posteriori error estimates for the DPG method. They are not optimal though as they do not take into account the fact that $\psi = 0$ and the Galerkin orthogonality of residual $l - R_V\psi_h - Bu_h$ to $\psi_h$. Using the first equation in (3.13), we obtain:

$$R_V(\underbrace{\psi}_{=0} - \psi_h) + B(u - u_h) = l - R_V\psi_h - Bu_h,$$

which, by the Galerkin orthogonality, implies:

$$\|B(u - u_h)\|_{V'}^2 = \|l - R_V\psi_h - Bu_h\|_{V'}^2 + \underbrace{\|\psi_h\|_V^2}_{=\|B'\psi_h\|_{U'}^2}. \tag{3.24}$$

This inequality, stronger than estimate (3.21), has served as a starting point for the a-posteriori error estimate presented in [3].

For the DPG* method, $l = 0$ implies $R_V \psi = -Bu$. Estimates (3.20) and (3.22) collapse to a single estimate:
$$\|Bu\|_{V'} = \|\psi\|_V \leq \|g\|_{U'}.$$
In fact, recall that $\psi_0 = 0$ so,
$$\|Bu\|_{V'} = \|\psi^\perp\|_V = \|g\|_{U'}.$$
However, this does not simplify error estimates (3.18),(3.21) and (3.23) as $R_V(\psi - \psi_h) \neq -B(u - u_h)$.

**Goal-oriented a-posteriori estimate.** Recall the classical reasoning for the standard (Petrov-) Galerkin method. Given a goal $g \in U'$, we want to estimate the error in goal:
$$g(u) - g(u_h) = g(u - u_h)$$
where $u$ and $u_h$ are the exact and approximate solutions of
$$\begin{cases} u \in U \\ b(u, \delta v) = l(\delta v) \quad \delta v \in V. \end{cases}$$
In order to relate the residual with the goal, we introduce the dual problem,
$$\begin{cases} v \in V \\ b(\delta u, v) = g(\delta u) \quad \delta u \in U. \end{cases}$$
We obtain then the standard identity:
$$\begin{aligned} g(u - u_h) &= b(u - u_h, v) && \text{(definition of dual solution } v\text{)} \\ &= b(u - u_h, v - v_h) && \text{(Galerkin orthogonality for the primal problem)} \end{aligned} \qquad (3.25)$$
where $v_h \in V_h$ is *an arbitrary element of* discrete space $V_h$. In particular, we can select for $v_h$ the discrete solution of the dual problem.

The sesquilinear form corresponding to the mixed problem (3.13),
$$\mathsf{b}(\mathsf{u}, \mathsf{v}) = \mathsf{b}((\psi, u), (\phi, w)) = (\psi, \phi)_V + b(u, \phi) + \overline{b(w, \psi)}$$
is Hermitian so, for the mixed formulation, the original and dual problems differ only in the load. Let $(\psi_h, u_h)$ be the approximate solution of (3.13) with $g = 0$, and $(\phi_h, w_h)$ the approximate solution of (3.13) with $l = 0$. We have,
$$\begin{aligned} g(u - u_h) &= (\psi - \psi_h, \phi - \phi_h)_V + b(u - u_h, \phi - \phi_h) + \overline{b(w - w_h, \psi - \psi_h)} \\ &= -(\psi_h, \phi - \phi_h)_V + b(u - u_h, \phi - \phi_h) - \overline{b(v - v_h, \psi_h)} && (\psi = 0) \\ &= b(u - u_h, \phi - \phi_h) - \overline{(\phi - \phi_h, \psi_h)_V + b(v - v_h, \psi_h)} \\ &= b(u - u_h, \phi - \phi_h) && \text{(Galerkin orthogonality)} \end{aligned}$$

We retain thus the standard Galerkin identity (3.25), even though our approximate solutions come from the mixed method. Another way to prove it is through the concept of the approximate optimal test functions. We have the Galerkin orthogonality,
$$b(u - u_h, v_h) = 0$$
for all *approximate optimal test functions*. The approximate solution of the dual problem satisfies
$$(\phi_h, \delta v_h) + b(w_h, \delta v_h) = 0$$
and, therefore, can be identified as the optimal test function corresponding to the approximate solution $w_h$. Consequently,
$$b(u - u_h, \phi) = b(u - u_h, \phi - \phi_h)$$
as claimed. For the actual results on goal-oriented estimation and adaptivity using DPG methods, see [10].

**Relation with weakly conforming least squares.** Let $A$ represent the operator corresponding to a first order system and $A'$ its $L^2$-adjoint. The weakly conforming least squares method for the adjoint problem seeks a minimizer of the least squares functional,
$$\frac{1}{2}\|A'\psi - g\|^2,$$
under the conformity constraint:
$$\langle \psi, \delta \hat{u} \rangle = 0 \quad \forall \delta \hat{u}$$
where, of course, operator $A'$ is now understood element-wise. This leads to a mixed problem:
$$\begin{cases} (A'\psi, A'\phi) + \langle \hat{u}, \phi \rangle &= (g, A'\phi) \quad \forall \phi \\ \langle \psi, \delta \hat{u} \rangle &= 0 \quad \quad \forall \delta \hat{u} \end{cases}$$
If we use the adjoint graph norm in the DPG* ultra-weak formulation, we get
$$\begin{cases} (A'\psi, A'\phi)_V + \alpha(\psi, \phi) + (u, A'\phi) + \langle \hat{u}, \phi \rangle &= 0 \quad \quad \forall \phi \in V \\ (A'\psi, \delta u) &= (g, \delta u) \quad \forall \delta u \in U \\ \langle \psi, \delta \hat{u} \rangle &= 0 \quad \quad \forall \delta \hat{u} \end{cases}$$
or,
$$\begin{cases} (A'\psi, A'\phi)_V + \alpha(\psi, \phi) + \langle \hat{u}, \phi \rangle &= (-u, A'\phi) \quad \forall \phi \in V \\ (A'\psi, \delta u) &= (g, \delta u) \quad \quad \forall \delta u \in U \\ \langle \psi, \delta \hat{u} \rangle &= 0 \quad \quad \forall \delta \hat{u}. \end{cases}$$
Passing with $\alpha \to 0$ we see that, in the limit, there should be $A'\psi = g = -u$. Consequently, this particular DPG* formulation can be viewed as a regularization of the weakly conforming least-squares. Figure 5 presents an eye-ball norm comparison of DPG* and least squares solutions for the linear acoustics problem with $p = 3$ and $dp = 3$.

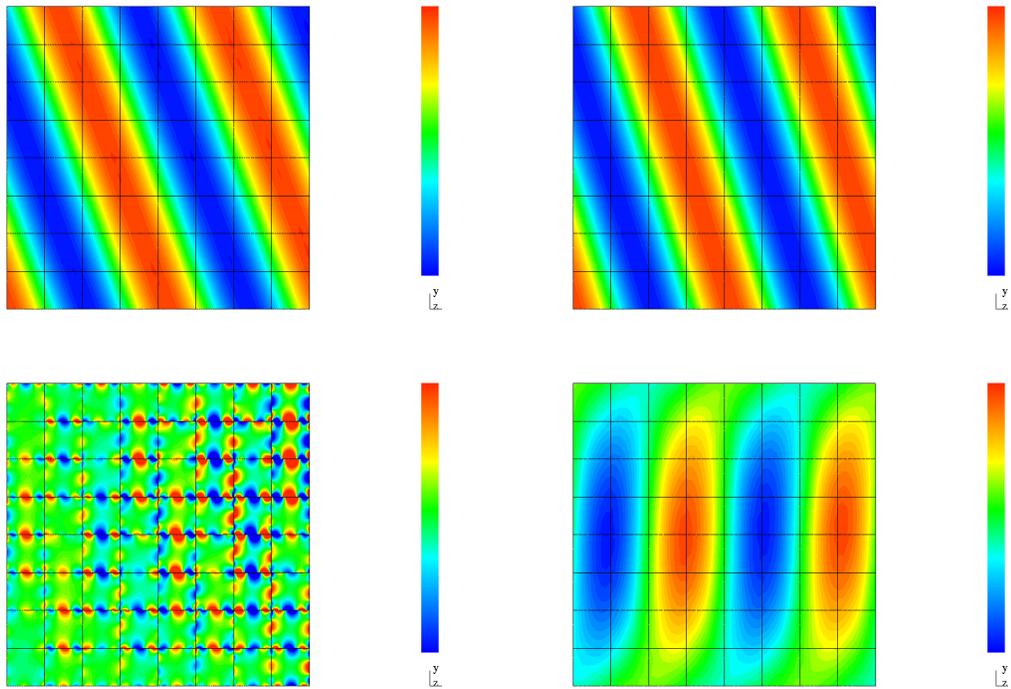

Figure 5: Plane wave example with 2 wavelengths. Comparison of DPG*(left) and weakly conforming least squares (right) solutions. Real part of pressure (top) and the corresponding pointwise error (bottom).

# 4 Conclusions

We have arrived at the DPG* method in the course of studying goal-oriented adaptivity and the need for a discretization of the dual problem $A'\psi = g$. Both methods stem from the discretization of the same mixed problem and share the same stiffness matrix. The existing stability and convergence theory for the DPG method automatically extends to the DPG* method. We have compared the two methods within a simple 2D linear acoustics problem and the ultraweak variational formulation. A particular motivation for choosing the problem is the (formal) relationship of DPG* methods with the weakly conforming least squares methods. Compared with DPG, the DPG* method seems to deliver slightly better results in $L^2$-norm but it is more sensitive to the round off error. It remains to be seen whether the DPG* method is only a "poor cousin" of the DPG method or it has a right to exist on its own.